\documentclass[12pt]{amsart} 
\usepackage{amsmath,amssymb,latexsym,amscd,graphicx}

\theoremstyle{definition}

\newtheorem{thm}{Theorem}
\newtheorem{cor}{Corollary}

\newtheorem{defi}{Definition}

\def\Z{\raise.5ex\hbox{$\chi$}}

\def\cs{{\mathcal S}}

\def\cl{\mbox{\textbf{cl}}}
\def\ext{\mbox{\textbf{ext}}}
\def\int{\mbox{\textbf{int}}}

\def\proof{\noindent{\textit{Proof.} \;}}

\begin{document}

\title[continuity and separation]{Continuity and separation for pointwise-symmetric isotonic closure functions}

\author{John M. Harris}
\address{Department of Mathematics\\
The University of Southern Mississippi\\
Long Beach, Mississippi}
\email{john.m.harris@usm.edu}

\begin{abstract}  In this paper, we show that a pointwise-symmetric isotonic closure function is uniquely determined by the pairs of sets it separates.  We then show that when the closure function of the domain is isotonic and the closure function of the codomain is isotonic and pointwise-symmetric, functions which separate only those pairs of sets which are already separated are continuous, generalizing a result in \cite{HL} which is, in turn, a generalization of a result in \cite{L}.
\end{abstract}

\maketitle

A \emph{generalized closure space} is a pair $(X, \cl)$ consisting of a set $X$ and a \emph{closure} function $\cl$, a function from the power set of $X$ to itself.  The \emph{closure} of a subset $A$ of $X$, denoted $\cl(A)$, is the image of $A$ under $\cl$.  The \emph{exterior} of $A$ is $\ext(A) = X - \cl(A)$, and the \emph{interior} of $A$ is $\int(A) = X - \cl(X - A)$.
\medskip

We say that $A$ is \emph{closed} if $A = \cl(A)$, $A$ is \emph{open} if $A = \int(A)$, and $N$ is a \emph{neighborhood} of $x$ if $x \in \int(N)$.
\medskip

In this paper, we use terminology for generalized closure spaces found in \cite{SS}, namely, we say that a closure function $\cl$ defined on $X$ is 
\medskip

\begin{enumerate}

\item \emph{grounded}, if $\cl(\emptyset) = \emptyset$.
\medskip

\item \emph{isotonic}, if $\cl(A) \subseteq \cl(B)$ whenever $A \subseteq B$.  
\medskip

\item \emph{enlarging}, if $A \subseteq \cl(A)$, for each subset $A$ of $X$.  
\medskip

\item \emph{idempotent}, if $\cl(A) = \cl(\cl(A))$, for each subset $A$ of $X$.
\medskip

\item \emph{sub-linear}, if $\cl(A \cup B) \subseteq \cl(A) \cup \cl(B)$, for all $A, B \subseteq X$.
  
\end{enumerate}
\medskip

Additionally, we define a kind of separation useful for the task at hand:

\begin{defi} Subsets $A$ and $B$ of $X$ are said to be \emph{closure-separated} in a generalized closure space $(X, \cl)$ (or simply, \emph{$\cl$-separated}) if $A \cap \cl(B) = \emptyset$ and $\cl(A) \cap B = \emptyset$, or, equivalently, if $A \subseteq \ext(B)$ and $B \subseteq \ext(A)$.  
\end{defi}

Before we begin, we also define a relatively weak separation axiom:

\begin{defi}
\emph{Exterior points are closure-separated} in a generalized closure space $(X, \cl)$ if, for each $A \subseteq X$ and for each $x \in \ext(A)$, $\{x\}$ and $A$ are $\cl$-separated.
\end{defi}

\begin{thm}
\label{clthm}
Let $(X, \cl)$ be a generalized closure space in which exterior points are $\cl$-separated, and let $\cs$ be the pairs of $\cl$-separated sets in $X$.  Then, for each subset $A$ of $X$, the closure of $A$ is $$\cl(A) = \{ x \in X : \{ \{x\}, A\} \notin \cs\}.$$
\end{thm}
\proof In any generalized closure space, $\cl(A) \subseteq \{ x \in X : \{ \{x\}, A\} \notin \cs\}$:  suppose that $y \notin \{ x \in X : \{ \{x\}, A\} \notin \cs\}$, that is, $\{ \{y\}, A\} \in \cs$.  Then $\{y\} \cap \cl(A) = \emptyset$, and so, $y \notin \cl(A)$.
\medskip

Suppose now that $y \notin \cl(A)$.  By hypothesis, $\{ \{y\}, A\} \in \cs$, and hence, $y \notin \{ x \in X : \{ \{x\}, A\} \notin \cs\}$.  \qed
\bigskip

\begin{defi}
A closure function $\cl$ defined on a set $X$ is \emph{pointwise-symmetric} when, for all $x, y \in X$, if $x \in \cl(\{y\})$, then $y \in \cl(\{x\})$.

A generalized closure space $(X, \cl)$ is $R_0$ when, for all $x, y \in X$, if $x$ is in each neighborhood of $y$, then $y$ is in each neighborhood of $x$. 
\end{defi}

In \cite{SS}, a generalized closure space with a pointwise-symmetric closure function is said to be $(R0c)$, while that which we denote by $R_0$ is there denoted $(R0)$.  Note that both conditions hold whenever exterior points are closure-separated:

\begin{cor}
Let $(X, \cl)$ be a generalized closure space in which exterior points are $\cl$-separated.  Then $\cl$ is pointwise-symmetric and $(X, \cl)$ is $R_0$.
\end{cor}
\proof Suppose that exterior points are $\cl$-separated in $X$.  If $x \in \cl(\{y\})$, then $\{x\}$ and $\{y\}$ are not $\cl$-separated, and hence, $y \in \cl(\{x\})$.  Hence, $\cl$ is pointwise-symmetric. 
\medskip

Suppose that $x$ belongs to every neighborhood of $y$, that is, $x \in M$ whenever $y \in \int(M)$.  Letting $A = X - M$ and rewriting contrapositively, $y \in \cl(A)$ whenever $x \in A$.  
\medskip

Suppose $x \in \int(N)$.  $x \notin \cl(X - N)$, so $x$ is $\cl$-separated from $X - N$.  Hence, $\cl(\{x\}) \subseteq N$.
\medskip

$x \in \{x\}$, so $y \in \cl(\{x\}) \subseteq N$.  Hence, $(X, \cl)$ is $R_0$.  \qed
\bigskip

While these three axioms are not equivalent in general, they are equivalent when the closure function is isotonic:

\begin{thm} Let $(X, \cl)$ be a generalized closure space with $\cl$ isotonic.  Then the following are equivalent:
\medskip

\begin{enumerate}
\item exterior points are $\cl$-separated.
\medskip

\item $\cl$ is pointwise-symmetric.
\medskip

\item $(X, \cl)$ is $R_0$.
\end{enumerate}
\end{thm}
\proof  Suppose that (2) is true.  Let $A \subseteq X$, and suppose $x \in \ext(A)$.  Then, as $\cl$ is isotonic, for each $y \in A$, $x \notin \cl(\{y\})$, and hence, $y \notin \cl(\{x\})$.  Hence, $\cl(\{x\}) \cap A = \emptyset$.  Hence, (2) implies (1), and by the previous corollary, (1) implies (2).
\medskip

Suppose now that (2) is true and let $x, y \in X$ such that $x$ is in every neighborhood of $y$, that is, $x \in N$ whenever $y \in \int(N)$.  Then $y \in \cl(A)$ whenever $x \in A$, and in particular, since $x \in \{x\}$, $y \in \cl(\{x\})$.  Hence, $x \in \cl(\{y\})$.  Thus, if $y \in B$, then $x \in \cl(\{y\}) \subseteq \cl(B)$, as $\cl$ is isotonic.  Hence, if $x \in \int(C)$, then $y \in C$, that is, $y$ is in every neighborhood of $x$.  Hence, (2) implies (3).
\medskip

Finally, suppose that $(X, \cl)$ is $R_0$ and suppose that $x \in \cl(\{y\})$.  Since $\cl$ is isotonic, $x \in \cl(B)$ whenever $y \in B$, or, equivalently, $y$ is in every neighborhood of $x$.  Since $(X, \cl)$ is $R_0$, $x \in N$ whenever $y \in \int(N)$.  Hence, $y \in \cl(A)$ whenever $x \in A$, and in particular, since $x \in \{x\}$, $y \in \cl(\{x\})$.  Hence, (3) implies (2).  \qed
\bigskip

Hence, when $\cl$ is isotonic and pointwise-symmetric, the collection of $\cl$-separated sets uniquely determines $\cl$.  In fact, such closure functions can be defined entirely in terms of the pairs of sets which they closure-separate:  

\begin{thm} Let $\cs$ be a set of unordered pairs of subsets of a set $X$ such that, for all $A, B, C \subseteq X$, 
\medskip

\begin{enumerate}
	\item if $A \subseteq B$ and $\{ B, C \} \in \cs$, then $\{ A, C \} \in \cs$, and 
	\medskip

	\item if $\{ \{x\}, B\} \in \cs$ for each $x \in A$ and $\{ \{y\}, A\} \in \cs$ for each $y \in B$, then $\{A, B\} \in \cs$.
\end{enumerate}
\medskip

Then there is a unique pointwise-symmetric isotonic closure function $\cl$ on $X$ which closure-separates the elements of $\cs$. 
\end{thm} 

\proof Define $\cl$ by $\cl(A) = \{ x \in X : \{ \{x\}, A \} \notin \cs\}$, for every $A \subseteq X$.  If $A \subseteq B \subseteq X$ and $x \in \cl(A)$, then $\{ \{x\}, A \} \notin \cs$.  Hence, $\{ \{x\}, B \} \notin \cs$, that is, $x \in \cl(B)$.  Hence, $\cl$ is isotonic.  Also, $x \in \cl(\{y\})$ iff $\{ \{x\}, \{y\} \} \notin \cs$ iff $y \in \cl(\{x\})$, and thus, $\cl$ is pointwise-symmetric.
\medskip

Suppose that $\{A, B\} \in \cs$.  Then $A \cap \cl(B) = A \cap \{ x \in X : \{ \{ x \}, B\} \notin \cs \} = \{ x \in A : \{ \{ x \}, B\} \notin \cs \} = \emptyset$.  Similarly, $\cl(A) \cap B = \emptyset$.  Hence, if $\{A, B\} \in \cs$, then $A$ and $B$ are $\cl$-separated.
\medskip

Now suppose that $A$ and $B$ are $\cl$-separated.  Then $\{ x \in A : \{ \{ x \}, B\} \notin \cs \} = A \cap \cl(B) = \emptyset$ and $\{ x \in A : \{ \{ x \}, B\} \notin \cs \} = \cl(A) \cap B = \emptyset$.  Hence, $\{ \{x\}, B\} \in \cs$ for each $x \in A$ and $\{ \{y\}, A\} \in \cs$ for each $y \in B$, and thus, $\{A, B\} \in \cs$.  \qed
\bigskip

Furthermore, many properties of closure functions can be expressed in terms of the sets they separate:

\begin{thm} Let $\cs$ be the pairs of $\cl$-separated sets of a generalized closure space $(X, \cl)$ in which exterior points are closure-separated.  Then $\cl$ is 
\medskip

\begin{enumerate}
\item grounded if, and only if, for all $x \in X$, $\{ \{x\}, \emptyset \} \in \cs$. 
\medskip

\item enlarging if, and only if, for all $\{ A, B \} \in \cs$, $A$ and $B$ are disjoint.
\medskip

\item sub-linear if, and only if, $\{A, B \cup C \} \in \cs$ whenever $\{A, B\} \in \cs$ and $\{A, C\} \in \cs$.
\end{enumerate}

Moreover, if $\cl$ is enlarging and for all $A, B \subseteq X$, $\{ \{x\}, A\} \notin \cs$ whenever $\{ \{x\}, B \} \notin \cs$ and $\{ \{y\}, A \} \notin \cs$, for each $y \in B$, then $\cl$ is idempotent.  Also, if $\cl$ is isotonic and idempotent, then $\{ \{x\}, A\} \notin \cs$ whenever $\{ \{x\}, B \} \notin \cs$ and $\{ \{y\}, A \} \notin \cs$, for each $y \in B$.
\end{thm}
\proof  Recall that by Theorem \ref{clthm}, $\cl(A) = \{ x \in X : \{ \{x\}, A \} \notin \cs\}$, for every $A \subseteq X$.  
\medskip

Suppose that for all $x \in X$, $\{ \{x\}, \emptyset \} \in \cs$.  Then $\cl(\emptyset) = \{ x \in X : \{ \{x\}, \emptyset \} \notin \cs \} = \emptyset$.  Hence, $\cl$ is grounded.  

Conversely, if $\emptyset = \cl(\emptyset) = \{ x \in X : \{ \{x\}, \emptyset \} \notin \cs \}$, then $\{ \{x\}, \emptyset \} \in \cs$, for all $x \in X$.
\medskip

Suppose that for all $\{ A, B \} \in \cs$, $A$ and $B$ are disjoint. Since  $\{\{a\}, A\} \notin \cs$ if $a \in A$, $A \subseteq \cl(A)$, for each $A \subseteq X$.  Hence, $\cl$ is enlarging.

Conversely, suppose that $\cl$ is enlarging and $\{ A, B \} \in \cs$.  Then $A \cap B \subseteq \cl(A) \cap B = \emptyset$.
\medskip

Suppose that $\{A, B \cup C \} \in \cs$ whenever $\{A, B\}, \{A, C\} \in \cs$.  Let $x \in X$ and $B, C \subseteq X$ such that $\{ \{x\}, B \cup C \} \notin \cs$.  Then $\{ \{x\}, B\} \notin \cs$ or $\{ \{x\}, C\} \notin \cs$.  
Hence, $\cl(B \cup C) \subseteq \cl(B) \cup \cl(C)$, and therefore, $\cl$ is sub-linear.

Conversely, suppose that $\cl$ is sub-linear, and let $\{A, B\}, \{A, C\} \in \cs$.  Then $\cl(A \cup B) \cap C \subseteq (\cl(A) \cup \cl(B)) \cap C = (\cl(A) \cap C) \cup (\cl(B) \cap C) = \emptyset$, and $(A \cup B) \cap \cl(C) = (A \cap \cl(C)) \cup (B \cap \cl(C)) = \emptyset$.
\medskip

Suppose that $\cl$ is enlarging, and suppose that $\{ \{x\}, A\} \notin \cs$ whenever $\{ \{x\}, B \} \notin \cs$ and $\{ \{y\}, A \} \notin \cs$, for each $y \in B$.  Then $\cl(\cl(A)) \subseteq \cl(A)$:  if $x \in \cl(\cl(A))$, then $\{ \{x\}, \cl(A)\} \notin \cs$.  $\{ \{y\}, A \} \notin \cs$, for each $y \in \cl(A)$;  hence, $\{ \{x\}, A\} \notin \cs$.  And since $\cl$ is enlarging, $\cl(A) \subseteq \cl(\cl(A))$.  Thus, $\cl(\cl(A)) = \cl(A)$, for each $A \subseteq X$.  

Finally, suppose that $\cl$ is isotonic and idempotent.  Let $x \in X$ and $A, B \subseteq X$ such that $\{ \{x\}, B \} \notin \cs$ and, for each $y \in B$, $\{ \{y\}, A \} \notin \cs$.  Then $x \in \cl(B)$ and for each $y \in B$, $y \in \cl(A)$, that is, $B \subseteq \cl(A)$.  Hence, $x \in \cl(B) \subseteq \cl(\cl(A)) = \cl(A)$.  \qed
\bigskip

\begin{defi} If $(X, \cl_X)$ and $(Y, \cl_Y)$ are generalized closure spaces, then a function $f:X \longrightarrow Y$ is said to be 
\medskip

\begin{enumerate}

\item \emph{closure-preserving}, if $f(\cl_X(A)) \subseteq \cl_Y(f(A))$, for each $A \subseteq X$. 
\medskip

\item \emph{continuous}, if $\cl_X(f^{-1}(B)) \subseteq f^{-1}(\cl_Y(B))$, for each $B \subseteq Y$. 
 
\end{enumerate}
\end{defi}

In general, neither condition implies the other.  However, we easily obtain the following result:

\begin{thm}
Let $(X, \cl_X)$ and $(Y, \cl_Y)$ be generalized closure spaces, and let $f:X \longrightarrow Y$.
\medskip

\begin{enumerate}
\item If $f$ is closure-preserving and $\cl_Y$ is isotonic, then $f$ is continuous.
\medskip

\item If $f$ is continuous and $\cl_X$ is isotonic, then $f$ is closure-preserving.

\end{enumerate}
\end{thm}
\proof Suppose that $f$ is closure-preserving and $\cl_Y$ is isotonic.  Let $B \subseteq Y$.  $f(\cl_X(f^{-1}(B))) \subseteq \cl_Y(f(f^{-1}(B))) \subseteq \cl_Y(B)$, and hence, $\cl_X(f^{-1}(B)) \subseteq f^{-1}(f(\cl_X(f^{-1}(B)))) \subseteq f^{-1}(\cl_Y(B))$.
\medskip

Suppose that $f$ is continuous and $\cl_X$ is isotonic.  Let $A \subseteq X$.  $\cl_X(A) \subseteq \cl_X(f^{-1}(f(A))) \subseteq f^{-1}(\cl_Y(f(A)))$, and hence, $f(\cl_X(A)) \subseteq f(f^{-1}(\cl_Y(f(A)))) \subseteq \cl_Y(f(A))$. \qed
\bigskip

\begin{defi} Let $(X, \cl_X)$ and $(Y, \cl_Y)$ be generalized closure spaces, and let $f:X \longrightarrow Y$.  If, for all $A, B \subseteq X$,  $f(A)$ and $f(B)$ are not $\cl_Y$-separated whenever $A$ and $B$ are not $\cl_X$-separated, then we say that $f$ is \emph{nonseparating}.
\end{defi}

Note that $f$ is nonseparating if and only if $A$ and $B$ are $\cl_X$-separated whenever $f(A)$ and $f(B)$ are $\cl_Y$-separated.

\begin{thm} Let $(X, \cl_X)$ and $(Y, \cl_Y)$ be generalized closure spaces, and let $f:X \longrightarrow Y$.  
\medskip

\begin{enumerate}
\item If $\cl_Y$ is isotonic and $f$ is nonseparating, then $f^{-1}(C)$ and $f^{-1}(D)$ are $\cl_X$-separated whenever $C$ and $D$ are $\cl_Y$-separated.\medskip

\item If $\cl_X$ is isotonic and $f^{-1}(C)$ and $f^{-1}(D)$ are $\cl_X$-separated whenever $C$ and $D$ are $\cl_Y$-separated, then $f$ is nonseparating.

\end{enumerate} 
\end{thm}
\proof Let $C$ and $D$ be $\cl_Y$-separated subsets, where $\cl_Y$ is isotonic.  Let $A = f^{-1}(C)$ and let $B = f^{-1}(D)$.  $f(A) \subseteq C$ and $f(B) \subseteq D$, and since $\cl_Y$ is isotonic, $f(A)$ and $f(B)$ are also $\cl_Y$-separated.  Hence, $A$ and $B$ are $\cl_X$-separated in $X$.
\medskip

Suppose that $\cl_X$ is isotonic and let $A, B \subseteq X$ such that $C = f(A)$ and $D = f(B)$ are $\cl_Y$-separated.  Then $f^{-1}(C)$ and $f^{-1}(D)$ are $\cl_X$-separated, and since $\cl_X$ is isotonic, $A \subseteq f^{-1}(f(A)) = f^{-1}(C)$ and $B \subseteq f^{-1}(f(B)) = f^{-1}(D)$ are $\cl_X$-separated as well.  \qed
\bigskip

\begin{thm}
\label{cp_implies_ns} Let $(X, \cl_X)$ and $(Y, \cl_Y)$ be generalized closure spaces, and let $f:X \longrightarrow Y$.  If $f$ is closure-preserving, then $f$ is nonseparating.
\end{thm}
\proof Suppose that $f$ is closure-preserving and that $A, B \subseteq X$ are not $\cl_X$-separated.  Suppose that $\cl_X(A) \cap B \neq \emptyset$.  Then $\emptyset \neq f(\cl_X(A) \cap B) \subseteq f(\cl_X(A)) \cap f(B) \subseteq \cl_Y(f(A)) \cap f(B)$.  Similarly, if $A \cap \cl_X(B) \neq \emptyset$, then $f(A) \cap \cl_Y(B) \neq \emptyset$.  Hence, $f(A)$ and $f(B)$ are not $\cl_Y$-separated.  \qed
\bigskip

\begin{cor} Let $(X, \cl_X)$ and $(Y, \cl_Y)$ be generalized closure spaces with $\cl_X$ isotonic, and let $f:X \longrightarrow Y$.  If $f$ is continuous, then $f$ is nonseparating.
\end{cor}
\proof If $f$ is continuous and $\cl_X$ is isotonic, then $f$ is closure-preserving.  Hence, by the previous result, $f$ is nonseparating. \qed
\bigskip

\begin{thm}
Let $(X, \cl_X)$ and $(Y, \cl_Y)$ be generalized closure spaces with exterior points $\cl_Y$-separated in $Y$, and let $f:X \longrightarrow Y$.  Then $f$ is nonseparating iff $f$ is closure-preserving.
\end{thm}
\proof By Theorem \ref{cp_implies_ns}, if $f$ is closure-preserving, then $f$ is nonseparating.  Suppose that $f$ is nonseparating, and let $A \subseteq X$.  If $\cl_X(A) = \emptyset$, then $f(\cl_X(A)) = \emptyset \subseteq \cl_Y(f(A))$. 
\medskip

Suppose $\cl_X(A) \neq \emptyset$.  Let $\cs_X$ and $\cs_Y$ denote the pairs of $\cl_X$-separated subsets of $X$ and the pairs of $\cl_Y$-separated subsets of $Y$, respectively.  Let $y \in f(\cl_X(A))$, and let $x \in \cl_X(A) \cap f^{-1}(\{y\})$.  Since $x \in \cl_X(A)$, $\{ \{x\}, A\} \notin \cs_X$, and since $f$ is nonseparating, $\{ \{y\}, f(A)\} \notin \cs_Y$.  Since exterior points are $\cl_Y$-separated, $y \in \cl_Y(f(A))$.  Thus, $f(\cl_X(A)) \subseteq \cl_Y(f(A))$, for each $A \subseteq X$.  \qed
\bigskip

\begin{cor}
Let $(X, \cl_X)$ and $(Y, \cl_Y)$ be generalized closure spaces with isotonic closure functions and with $\cl_Y$ pointwise-symmetric, and let $f:X \longrightarrow Y$.  Then $f$ is nonseparating iff $f$ is continuous.
\end{cor}

\proof Since $\cl_Y$ is isotonic, exterior points are closure-separated in $(Y, \cl_Y)$.  Since both closure functions are isotonic, $f$ is closure-preserving iff $f$ is continuous.  Hence, we can apply the previous theorem.  \qed

\end{document}